\documentclass[12pt]{article}
\usepackage{amsfonts,latexsym,amssymb}
 \newcommand{\beqn}{\begin{eqnarray}}
 \newcommand{\eeqn}{\end{eqnarray}}
 \newcommand{\be}{\begin{equation}}
 \newcommand{\ee}{\end{equation}}
 \newcommand{\ba}{\begin{array}}
 \newcommand{\ea}{\end{array}}
 \newcommand{\pa}{\partial}
 
 \newcommand{\ci}{\cite}
 \newcommand{\la}{\label}

\newcommand{\ga}{\gamma}

\newcommand{\na}{\nabla}

 \newcommand{\de}{\delta}

 \newcommand{\De}{\Delta}

\def\R{{\rm I\kern-.1567em R}}
\def\M{{\rm I\kern-.1567em M}}

\def\div {{\rm div}}

\def\ess{{\rm ess}}

 \newtheorem{theorem}{Theorem}[section]
 
 \newtheorem{definition}[theorem]{Definition}
 
 \newtheorem{lemma}[theorem]{Lemma}

 \newtheorem{cor}[theorem]{Corollary}
 \newtheorem{pro}[theorem]{Proposition}

\begin{document}

\begin{center} {\bf Regularity for Suitable Weak Solutions to the
Navier-Stokes Equations in Critical Morrey Spaces}\\
  \vspace{1cm}
 {\large
  G. Seregin}

 \end{center}
 \vspace{1cm}
 \noindent
 {\bf Abstract }
A  class of sufficient conditions of local regularity for suitable
weak
 solutions to the nonstationary three-dimensional  Navier-Stokes equations
 are discussed. The corresponding results are formulated
  in terms of functionals which are invariant
 with respect to the Navier-Stokes equations scaling. The famous
 Caffarelli-Kohn-Nirenberg condition is contained in that class as a particular case.

 \vspace {1cm}

\noindent {\bf 1991 Mathematical subject classification (Amer.
Math. Soc.)}: 35K, 76D.


\setcounter{equation}{0}
\section{Introduction and Main Result }

In the present paper, we address the problem of smoothness of a
certain class of weak solutions to the nonstationary
three-dimensional Navier-Stokes equations (NSE's)
\begin{equation}\label{11}
    \pa_t v+v\cdot \na v-\De v=-\na p,\qquad \div v=0.
\end{equation}
In our setting, they are considered in the unit space-time cylinder
$Q=B\times ]-1,0[$, where $B$ is the unit ball of $\mathbb R^3$
centered at the  origin. As usual, $v$ and $p$ stand for the
velocity field  and for the pressure field, respectively. We ask the
following question. What are minimal conditions which ensure
regularity of the velocity field $v$ at the space-time origin
$z=(x,t)=(0,0)=0$? Our definition of regularity of $v$ at the point
$z=0$ means that there exists a number $r\in ]0,1]$ such that $v$ is
a H\"older continuous function in the completion of the space-time
cylinder $Q(r)=B(r)\times ]-r^2,0[$. Here, $B(r)$ is the ball of
radius $r$ centered at the origin so that $B=B(1)$. This definition
is due to O. Ladyzhenskaya and the author, see \ci{LS}, and slightly
differs from the most popular one by L. Caffarelli, R.-V. Kohn, and
L. Nirenberg given in their celebrated paper \ci{CKN}.

Our interest to the above question is motivated by the important
observation of J. Leray made in his remarkable paper \ci{Le}. J.
Leray proved the uniqueness of regular solutions in the class of
turbulent solutions which are also called weak Leray-Hopf
solutions. For exact definitions, we refer the reader to the paper
\ci{Le} and, for example, to papers \ci{LS} and \ci{S7}.

We restrict ourselves to the analysis of the so-called suitable
weak solutions also introduced in  \ci{CKN}.  Here, we follow the
definition of suitable weak solution in F.-H. Lin's reduction, see
\ci{Li} and also \ci{LS} for discussions on other definitions.
\begin{definition}\la{1d1} We say that the pair $v$ and $p$ is a suitable
weak solution to the NSE's in $Q$ if the following conditions are
fulfilled:
\begin{equation}\label{12}
    v\in L_{2,\infty}(Q)\cap W^{1,0}_2(Q),
\qquad p\in L_\frac 32(Q);
\end{equation}
\begin{equation}\label{13}
    \mbox{\it{the NSE's hold in $Q$ in the sense of distributions}};
    \end{equation}
for a.a. $t\in ]-1,0[$, the local energy inequality
$$\int\limits_B\varphi(x,t)|v(x,t)|^2dx+2\int\limits_{-1}^t\int\limits_B
\varphi|\na v|^2dxdt'\leq\int\limits_{-1}^t\int\limits_B\Big(|v|^2
(\De \varphi+\pa_t \varphi)$$
\begin{equation}\label{14}
    +v\cdot \varphi(|v|^2+2p)\Big)dxdt'
\end{equation}
holds for all non-negative functions $\varphi\in
C^\infty_0(\mathbb R^3\times\mathbb R^1)$ vanishing in a
neighborhood of the parabolic boundary $\pa'Q$ of the cylinder
$Q$.
\end{definition}
Here, the following notion has been used:
$$L_{m,n}(Q)=L_n(-1,0;L_m(B)),\qquad
W^{1,0}_2(Q)=L_2(-1,0;W^1_2(B)),$$ and $L_m(B)$ and $W^1_m(B)$ are
the usual Lebesgue and Sobolev spaces, respectively.

Let us comment Definition \ref{1d1} briefly. The most essential
part of it is  local energy inequality (\ref{14}). Unfortunately,
we do not know whether or not any weak Leray-Hopf solution to the
initial boundary value problem for the NSE's satisfies the local
energy inequality but at least one of them does so. Moreover, in
the case of local (in time) strong solvability of that
problem,
 its solution belongs to the class of
suitable weak solutions at least up to the moment of time when the
first singularity occurs. Our assumption on the pressure field
$p$, see (\ref{12}), is motivated by the linear theory. In this
sense, it is valid for any weak Leray-Hopf solution if the data of
the initial boundary value problem are not too bad. Moreover, we
may vary classes for the pressure filed. The space $L_\frac 32$
seems to be the most convenient  for treating. For more details,
we refer the reader to the paper \ci{LS}.

It is known that system (\ref{11}) is invariant with respect to
the scaling
$$v^\lambda(x,t)=\lambda v(\lambda x,\lambda^2 t),\qquad
p^\lambda(x,t)=\lambda^2 p(\lambda x,\lambda^2 t).$$
 We call this scaling the natural one.

In the local regularity theory, functionals being invariant to the
natural scaling play a very important role. Here, it the list of
some of them:
$$A(r)=\mbox{ess}\sup\limits_{-r^2<t<0}\frac 1r\int\limits_{B(r)}
|v(x,t)|^2dx, \qquad E(r)=\frac 1r\int\limits_{Q(r)}|\na v|^2dz,$$
$$C(r)=\frac 1{r^2}\int\limits_{Q(r)}| v|^3dz,\qquad
H(r)=\frac 1{r^3}\int\limits_{Q(r)}| v|^2dz,$$
$$D_0(r)=\frac 1{r^2}\int\limits_{Q(r)}|p-[p]_{B(r)}|^\frac 32dz,$$
where
$$[p]_{B(r)}=\frac 1{|B(r)|}\int\limits_{B(r)}p(x,t)dx.$$
These functionals may be used to produce norms of special Morrey
classes. We call them critical Morrey spaces. The list of scaling
invariant functionals can be extended. For example, the norm
$M_{s,l}(r)=\|v\|_{L_{s,l}(Q(r))}$ is invariant with respect to
the natural scaling if $3/s+2/l=1$.

The most of results in the local regularity theory is formulated
with the help of those functionals and have the form of the
so-called $\varepsilon$-regularity conditions. A typical
$\varepsilon$-regularity condition reads:  if the norm of the
velocity field $v$ in critical Morrey's space is small enough,
then the space-time origin is a regular point of $v$. This is not
completely rigorous statement but reflects the spirit of the
$\varepsilon$-regularity theory quite well. For example, the
famous Caffarelli-Kohn-Nirenberg condition can be formulated as
follows.
\begin{theorem}\la{1t2} Let the pair $v$ and $p$ be
a suitable weak solution to the NSE's in $Q$. There is a universal
constant $\varepsilon_0$ such that if
\be\la{15}\sup\limits_{0<r\leq 1}E(r)<\varepsilon_0,\ee then $z=0$
is a regular point of $v$.\end{theorem}
Another important example is the local version of the
Ladyzheskaya-Prodi-Serrin condition (LPS-condition) proved by
Serrin and Struwe in papers \ci{Se} and \ci{Str0}. In our
situation, it reads
\begin{theorem}\la{1t3} Let the pair $v$ and $p$ be
a suitable weak solution to the NSE's in $Q$. Assume that numbers
$s$ and $l$ satisfy conditions $$\frac 3s+\frac 2l=1, \qquad 3\leq
s\leq +\infty,\quad 2\leq l\leq +\infty.$$ There is a positive
constant $\overline{\varepsilon}_0$ depending on $s$ and $l$ only
such that if $$M_{s,l}(1)=\sup\limits_{0<r\leq
1}M_{r,l}(r)<\overline{\varepsilon}_0,$$ then $z=0$ is a regular
point of $v$.\end{theorem}
Other  examples can be founded, for instance, in the paper \ci{S7}.

If we believe that suitable weak solutions are smooth, then it would
be natural to get rid of smallness of scaling invariant functionals
and show that their boundedness is sufficient for regularity. In a
number of cases, it is obvious. For example, in the case of
LPS-condition, if $s>3$,  boundedness of $M_{s,l}(1)$, together with
the absolute continuity of Lebesgue's integral and the natural
scaling, allows us to assume that $M_{s,l}(1)$ is small as we wish.
In the marginal case $s=3$ and $l=+\infty$, the above mentioned
reduction is much more subtle and based on backward uniqueness
results for the heat operator with variable lower order terms, see
\ci {SS1} and \ci{ESS4}. For functionals $A(r)$, $C(r)$, and $E(r)$,
this is an open problem, i.e., it is unknown whether or not their
boundedness implies regularity.


To formulate our main result, let us introduce the additional
notation
$$G=\inf\{\limsup_{r\to 0}E(r),\limsup_{r\to 0}C(r),\limsup_{r\to 0}A(r)\},$$
$$g=\inf\{\liminf_{r\to 0}E(r),\liminf_{r\to 0}C(r),
\liminf_{r\to 0}A(r),\liminf_{r\to 0}H(r),\liminf_{r\to
0}D_0(r)\}.$$ The main result of the paper is
\begin{theorem}\la{1t4}
Assume that the pair $v$ and $p$ is a suitable weak solution to
the NSE's in $Q$. For any $M>0$,  there exists a positive number
$\varepsilon(M)$ with the following property. If $G<M$ and
$g<\varepsilon(M)$, then $z=0$ is a regular point of
$v$.\end{theorem} Let us discuss simple consequences of Theorem
\ref{1t4}. In particular, one has the following generalization of
Theorem  \ref{1t2}.
\begin{cor}\la{1c5} Assume that the pair $v$ and $p$ is a suitable weak solution to
the NSE's in $Q$. Given $M>0$, let $\varepsilon(M)$ be the number
of Theorem \ref{1t4}. If, for some $M>0$, two conditions
$\limsup_{r\to 0}E(r)<M$ and $\liminf_{r\to 0}E(r)<\varepsilon(M)$
hold, then $z=0$ is a regular point of $v$.\end{cor} \noindent To
see that  Theorem \ref{1t2} can be deduced from Corollary
\ref{1c5}, it is sufficient to let $M=1$ and
$\varepsilon_0=\min\{1,\varepsilon(1)\}$. Another direct
consequence of Theorem \ref{1t4} can be stated as follows.
\begin{cor}\la{1c6}Assume that the pair $v$ and $p$ is a suitable weak
solution to the NSE's in $Q$.  If $G<+\infty$ and $g=0$, then
$z=0$ is a regular point of $v$.\end{cor} Certainly, there are
other versions of Theorem \ref{1t4}. For example, we have the
following  statement.
\begin{theorem}\la{1t7}Assume that the pair $v$ and $p$ is a suitable
weak solution to the NSE's in $Q$. For any $M>0$,  there exists a
positive number $\widehat{\varepsilon}(M)$ with the following
property. Let $G<M$ and
$$
\liminf_{r\to 0}E_3(r)<\widehat{\varepsilon}(M),$$ where
$$E_3(r)=\frac 1r\int\limits_{Q(r)}|v_{,3}|^2dz,\qquad
v_{,3}=\frac {\pa v}{\pa x_3}.$$ Then $z=0$ is a regular point of
$v$.\end{theorem}

  \noindent\textbf{Acknowledgement} The work was supported by the
Alexander von Humboldt Foundation, by  the RFFI grant
05-01-00941-a, and by the CRDF grant RU-M1-2596-ST-04.

\setcounter{equation}{0}
\section {Estimates of Suitable Weak Solutions to the NSE's}

In this section, we would like to present the list of the main
estimates of suitable weak solutions to the NSE's in $Q$. The
first of them is a consequence of  multiplicative inequalities and
has the form \be\la{21}C(r)\leq c\Big [\Big(\frac \varrho
r\Big)^3A^\frac 34(\varrho) E^\frac 34(\varrho)+\Big(\frac  r
\varrho \Big)^3 A^\frac 32(\varrho)\Big]\ee for all
$0<r\leq\varrho\leq 1$. Here and in what follows, we denote all
positive universal constants by $c$. A proof of (\ref{21}) is
given,  for example, in \ci{LS}.

There are two consequences of the local energy inequality:
\be\la{22}A(R/2)+E(R/2)\leq c\Big[C^\frac 23(R)+C(R) +C^\frac
13(R)D_0^\frac 23(R)\Big]\ee and \be\la{23}A(R/2)+E(R/2)\leq
c\Big[C^\frac 23(R) +C^\frac 13(R)D_0^\frac 23(R)+A^\frac
12(R)C^\frac 23(R)E^\frac 12(R)\Big]\ee for all $0<R\leq 1$.
Inequality (\ref{22}) follows from the local energy inequality
directly, a proof of (\ref{23}) can be found in \ci{LS}.

There are also three versions of the decay estimate for the
pressure: \be\la{24}D_0(r)\leq c\Big[\Big(\frac
r\varrho\Big)^\frac 52D_0(\varrho) +\Big(\frac \varrho
r\Big)^2C(\varrho)\Big]\ee or \be\la{25}D_0(r)\leq
c\Big[\Big(\frac r\varrho\Big)^\frac 52D_0(\varrho) +\Big(\frac
\varrho r\Big)^2A^\frac 12(\varrho)E(\varrho)\Big]\ee or
\be\la{26}D_0(r)\leq c\Big[\Big(\frac r\varrho\Big)^\frac
52D_0(\varrho) +\Big(\frac \varrho r\Big)^3A^\frac
34(\varrho)E^\frac 34(\varrho)\Big].\ee Inequality
(\ref{24})--(\ref{26}) are valid for all $0<r\leq\varrho\leq 1$.
Inequality (\ref{24}) is proved in \ci{S6} and inequalities
(\ref{25}) and (\ref{26})  are discussed in \ci{S8}.

The following lemma shows that if one of the quantities
$\sup_{m<r\leq 1}E(r)$, $\sup_{m<r\leq 1}C(r)$, or $\sup_{m<r\leq
1}A(r)$ is bounded, then all others are bounded too.

\begin{lemma}\la{2l1} Suppose that the pair $v$ and $p$ is a suitable
weak solution to the Navier-Stokes equations in $Q$. The following
estimates are valid:

a) Let \be\la{27}\sup\limits_{0<r\leq 1}E(r)=E_0<+\infty.\ee Then
there exists a positive constant $d$ depending only on $E_0$ such
that \be\la{28}A^\frac 32(r)+C(r)+D^2_0(r)\leq d(E_0)\Big(r^\frac
12 (A^\frac 32(1)+D_0^2(1))+1\Big)\ee for all $0<r\leq 1/4$;

b) Let \be\la{29}\sup\limits_{0<r\leq 1} C(r)=C_0<+\infty.\ee Then
\be\la{210}A(r)+D_0(r)+E(r)\leq c\Big(r^2D_0(1)+C_0+C_0^\frac
23\Big)\ee for all $0<r\leq 1/2$;

c) Let \be\la{211}\sup\limits_{0<r\leq 1}A(r)=A_0<+\infty.\ee Then
there exists a positive constant $e$ depending only on $A_0$ such
that \be\la{212}C^\frac 43(r)+D_0(r)+E(r)\leq
e(A_0)\Big(r^2(D_0(1)+E(1))+1\Big)\ee for all $0<r\leq 1/2$.
\end{lemma}
A proof Lemma \ref{2l1} is based upon estimates
(\ref{21})--(\ref{26}) and presented in \ci{S8}.

\setcounter{equation}{0}
\section {Proof of Theorem \ref{1t4}}

The key part of the proof of Theorem \ref{1t4} is

\begin{pro}\la{3p1} Let the pair $v$ and $p$ be a suitable weak
solution to the NSE's in $Q$.

For any $M>0$, there exists a positive number $\varepsilon_{10}=
\varepsilon_{10}(M)$ with the following property. If
\begin{equation}\label{31}
    \sup\limits_{0<r\leq 1} E(r)=E_0\leq M
\end{equation}
and
\begin{equation}\label{32}
    g_{r_*}=\min\{E(r_*),\,A(r_*),\,C(r_*),\,H(r_*),\,
    D_0(r_*)\}<\varepsilon_{10}(M)
\end{equation}
for some $r_*\in ]0,\min\{1/4,(A^\frac 32(1)+D_0^2(1))^{-2}\}[$,
then $z=0$ is a regular point of $v$.\end{pro} \textsc{Proof}
Assume that the statement of the proposition is false. Then there
exist a positive number $M$ and
 a sequence of suitable weak solutions $v^n$ and $p^n$ to
the NSE's in $Q$ such that, for any $n\in \mathbb N$, the
following two conditions hold
\begin{equation}\label{33}
    E(v^n,r)=\frac 1r\int\limits_{Q(r)}|\na v^n|^2dz\leq M
\end{equation}
for all $0<r\leq 1$
and
$$g_{r_n}(v^n,p^n)=$$\begin{equation}\label{34}
    =\min\{E(v^n,r_n),\,A(v^n,r_n),\,C(v^n,r_n),\,
    H(v^n,r_n),\,D_0(p^n,r_n)\}\leq \frac 1n
\end{equation}
for some
\begin{equation}\label{35}
    r_n\in ]0,\min\{1/4,(A^\frac 32(v^n,1)+D_0^2(p^n,1))^{-2}\}[,
\end{equation}
but $z=0$ is a singular point of $v^n$. Here, we have used the
notation
$$A(v^n,r)=\ess \sup\limits_{-r^2<t<0}\frac 1r\int\limits_{B(r)}
|v(x,t)|^2dx, \qquad C(v^n,r)=\frac 1{r^2}\int\limits_{Q(r)}
|v^n|^3dz,$$
$$H(v^n,r)=\frac 1{r^3}\int\limits_{Q(r)}
|v^n|^2dz,\qquad D_0(p^n,r)= \frac 1{r^2}\int\limits_{Q(r)} |p^n-
[p^n]_{B(r)}|^\frac 32dz.$$ On the other hand, since $z=0$ is a
singular point of $v^n$,
 there exists a universal positive number
$\varepsilon$ such that
\begin{equation}\label{36}
C(v^n,r)+D_0(p^n,r)\geq\varepsilon>0
\end{equation}
for all $0<r\leq 1$,
see, for example, \ci{LS}. We emphasize that (\ref{36}) is valid
for any natural number $n$.

By Lemma \ref{2l1} and by the properties of $r_n$, see (\ref{35}),
we find the  estimate
$$A^\frac 32(v^n,r)+C(v^n,r)+D_0^2(p^n,r)\leq$$\be\la{37}\leq d(M)\Big[
\Big(\frac r{r_n}\Big)^\frac 12 r_n^\frac 12(A^\frac 32(v^n,1)
+D_0^2(p^n,1))+1\Big]\leq d_0(M)\ee for all $r\in]0,r_n[$.

Now, let us scale our functions $v^n$ and $p^n$ so that
$$u^n(y,s)=r_nv^n(r_ny,r_n^2s),\qquad q^n(y,s)=r^2_np^n(r_ny,r_n^2s)
$$
for $(y,s)\in Q$. By the invariance of functionals and equations
with respect to the natural scaling,  we have:
\be\la{38}\mbox{\it{the pair}}\,\,u^n\,\,\mbox{ \it{and}}\,\,q^n
\,\,\mbox{ \it{is a suitable weak solution to the NSE's in}}\,\,
Q\ee for each $n\in \mathbb N$; \be\la{39}E(u^n,r)\leq M\ee for
all $0<r\leq 1$ and for each $n\in \mathbb N$;
\be\la{310}g_{r_n}(v^n,p^n)=g_1(u^n,q^n)\to 0\qquad as \quad n\to
+\infty;\ee \be\la{311}C(u^n,r)+D_0(q^n,r)\geq \varepsilon>0\ee
for all $0<r\leq 1$ and  for each $n\in \mathbb N$;
\be\la{312}A^\frac 32(u^n,r)+C(u^n,r)+D_0^2(q^n,r)\leq d_0(M)\ee
for all $0<r\leq 1$ and  for each $n\in \mathbb N$.

Now, let $n$ tend to $+\infty$. First of all, in order to pass to
the limit in non-linear terms,  a strong compactness is needed. To
this end, we estimate the weak derivative of $v$ in $t$ using the
equation in the standard way:
$$\int\limits_Q\pa_tu^n\cdot wdz=\int\limits_Q\Big(u^n\otimes
u^n:\na w -\na u^n:\na w+(q^n-[q^n]_B)\div w\Big)dz\leq$$
$$\leq \Big(\int\limits_Q|u^n|^3dz\Big)^\frac 23
\Big(\int\limits_Q|\na w|^3dz\Big)^\frac 13+\Big(\int\limits_Q|\na
u^n|^2dz \Big)^\frac 12 \Big(\int\limits_Q|\na w|^2dz\Big)^\frac
12+$$
$$ +\Big(\int\limits_Q|q^n-[q^n]|^\frac 32dz\Big)^\frac 23
\Big(\int\limits_Q|\div w|^3dz\Big)^\frac 13\leq$$
$$\leq c\Big(C^\frac 23(u^n,1)+D_0^\frac 23(q^n,1)+E^\frac 12(u^n,1)
\Big)\Big(\int\limits_Q|\na w|^3dz\Big)^\frac 13\leq$$
$$\leq d_1(M)\Big(\int\limits_Q|\na w|^3dz\Big)^\frac 13.$$
The latter estimate holds for any $w\in C^\infty_0(Q)$. By the
density arguments, it is valid for any $w\in
L_3(-1,0;{\stackrel{\circ}{W}}{^1_3}(B))$. So, we have
\be\la{313}\|\pa_t u^n\|_{L_{ 3/2}(-1,0;W^{-1}_{3/2}(B)
)}\leq d_1(M).\ee The final estimate comes from the known
multiplicative inequality and has the form
\be\la{314}\int\limits_Q |u^n|^\frac {10}3dz\leq d_2(M).\ee

Now, using known compactness arguments and selecting a subsequence
if necessary, we find
$$u^n\stackrel{\star}{\rightharpoondown}u\qquad in\quad L_{2,\infty}(Q),$$
\begin{equation}\label{315}
\na u^n\rightharpoondown \na u\qquad in \quad L_2(Q),
\end{equation}
$$u^n\rightarrow u\qquad in\quad L_3(Q),$$
$$q^n-[q^n]_B\rightharpoondown q\qquad in \quad L_\frac 32(Q).$$
Moreover, the pair $u$ and $q$ ia a suitable weak solution to the
NSE's in $Q$ and
\begin{equation}\label{316}
    A^\frac 32(u,r)+C(u,r)+D_0^2(q,r)\leq d_0(M)
\end{equation}
for all $0<r\leq 1$.

Now, let us see what follows from (\ref{310}). It is easy to
observe that there exists a subsequence $\{n_k\}_{k=1}^\infty$
such that either
\begin{equation}\label{317}
    E(u^{n_k},1)\rightarrow 0,
\end{equation}
or
\begin{equation}\label{318}
    C(u^{n_k},1)\rightarrow 0,
\end{equation}
or
\begin{equation}\label{319}
    A(u^{n_k},1)\rightarrow 0,
\end{equation}
or
\begin{equation}\label{320}
    H(u^{n_k},1)\rightarrow 0,
\end{equation}
or
\begin{equation}\label{321}
    D_0(q^{n_k},1)\rightarrow 0
\end{equation} as $k\rightarrow +\infty$.

Let us discuss each case separately, starting with (\ref{317}).
According to (\ref{317}), we have
\begin{equation}\label{322}
    \na u^{n_k}\rightarrow 0\qquad in \quad L_2(Q)
\end{equation}
and $\na u=0$ in $Q$. Therefore, $E(u,r)=0$ for all $r\in ]0,1]$
and, by Theorem \ref{1t2}, $z=0$ is a regular point of $u$ and, in
particular, there exists a number $0<r_1\leq 1$ such that
\begin{equation}\label{323}
    \sup\limits_{z\in \overline{Q}(r_1)}|u(z)|\leq d_3(M).
\end{equation}
Passing to the limit in (\ref{311}), we show
\be\la{324}C(u,r)+\limsup\limits_{k\to \infty}D_0(q^{n_k},r)\geq
\varepsilon\ee for all $0<r\leq 1$. Thanks to (\ref{323}), it
follows from (\ref{324}) that
\begin{equation}\label{325}
    cd^3_3(M)r^3+\limsup\limits_{k\to \infty}D_0(q^{n_k},r)\geq
\varepsilon
\end{equation}
for all $0<r\leq r_1$.

Functions $u^{n_k}$ and $q^{n_k}$ satisfy decay estimate
(\ref{25}). So, according to (\ref{312}), we have
$$D_0(q^{n_k},r)\leq c \Big[\Big(\frac r1\Big)^\frac 52
D_0(q^{n_k},1)+\Big(\frac 1r\Big)^2A^\frac
12(u^{n_k},1)E(u^{n_k},1)\Big]$$
$$\leq d_4(M)\Big[r^\frac 52+\Big(\frac 1r\Big)^2E(u^{n_k},1)\Big]$$
for all $0<r\leq 1$. It remains to take the limit as $k\to\infty$
and, by (\ref{322}), arrive at the inequality
$$\limsup\limits_{k\to \infty}D_0(q^{n_k},r)\leq d_4(M)r^\frac 52$$
for all $0<r\leq 1$. So, the latter estimate, together with
(\ref{325}), gives us the relation
$$cd^3_3(M)r^3+d_4(M)r^\frac 52\geq
\varepsilon $$ for all $0<r\leq r_1$. It is not true for
sufficiently small $r$.
So, possibility of (\ref{317}) is excluded.

Now, assume that (\ref{318}) takes place. Here, we are going to
use inequality (\ref{312}) and decay estimate (\ref{24}). They
lead us to the estimate
$$D_0(q^{n_k},r)\leq c \Big[ r^\frac 52
D_0(q^{n_k},1)+\frac 1{r^2}C(u^{n_k},1)\Big]$$
\begin{equation}\label{326}
\leq d_5(M)\Big[ r^\frac 52 +\frac 1{r^2}C(u^{n_k},1)\Big]
\end{equation}
for all $0<r\leq 1$. But it follows from (\ref{311}) that
\begin{equation}\label{327}
    \varepsilon\leq d_5(M)\Big[ r^\frac 52 +\frac 1{r^2}C(u^{n_k},1)\Big]
 +   C(u^{n_k},r),\qquad r\in ]0,1].
\end{equation}
Taking into account (\ref{318}) and passing to the limit in
(\ref{327}) as $k\to +\infty$,
 we show
$$\varepsilon\leq d_5(M) r^\frac 52, \qquad r\in ]0,1],$$
which is also wrong for sufficiently small $r$. So, case
(\ref{318}) is excluded as well.

Cases (\ref{319}) and (\ref{320}) can be reduced to the previous
one. Indeed, in both cases the limit equation is $u=0$. So,
$C(u^{n_k},r)\to C(u,r)=0$ as $k\to +\infty$ for all $0<r\leq 1$.
 Repeating estimates (\ref{326})
and (\ref{327}), we see that cases  (\ref{319}) and (\ref{320}) do
not occur either.

It remains to consider case (\ref{321}). Here, we have $q^{n_k}-
[q^{n_k}]\to 0$ in $L_\frac 32(Q)$ as $k\to +\infty$. Since
$$0\leq D_0(q^{n_k},r)\leq \frac c{r^2}D_0(q^{n_k},1)\to 0,$$
we find \be\la{328}C(u,r)\geq \varepsilon>0\ee for all $0<r\leq
1$. Let us describe the properties of the limit function $u$:
\be\la{329}u\in L_{2,\infty}(Q)\cap W^{1,0}_2(Q);\ee the function
$u$ satisfies the system of equations \be\la{303}\pa_t u+u\cdot
\na u-\De u=0,\qquad \div u=0\ee in $Q$ in the sense of
distributions;

\noindent for a.a. $t\in ]-1,0[$, the function $u$ satisfies the
local energy inequality
$$\int\limits_B\varphi(x,t)|u(x,t)|^2dx
+2\int\limits_{-1}^t\int\limits_B\varphi |\na u|^2dxdt'\leq$$
\be\la{331} \leq \int\limits_{-1}^t\int\limits_B\Big(|u|^2(\De
\varphi+\pa_t\varphi)+|u|^2u\cdot\na\varphi \Big)dxdt'\ee for all
non-negative functions $\varphi\in C_0^\infty(\mathbb R^3\times
\mathbb R)$ vanishing in a neighborhood of the parabolic boundary
of the cylinder $Q$.

As it is shown  in the Appendix, the function $u$, enjoying
properties (\ref{329})-(\ref{331}), is, in fact, smoother. More
precisely, $u$  is H\"older continuous say in $\overline{Q}(1/4)$
and, in particular,
$$\sup\limits_{z\in \overline{Q}(1/4)}|u(z)|\leq d_6(M)$$
and, from (\ref{328}), it follows that
$$cd^3_6(M)r^3\geq\varepsilon$$
for all $0<r\leq 1/4$, which is also not true. So, case
(\ref{321}) is excluded. Since our observations contradict with
(\ref{310}), we may conclude that the statement of Proposition
\ref{3p1} is valid. Proposition \ref{3p1} is proved.
\begin{pro}\la{3p2} Let the pair $v$ and $p$ be a suitable weak solution
to the NSE's in $Q$. If \be\la{332}\limsup\limits_{r\to 0}
E(r)<\frac 12m=M\ee and \be\la{333} g<\frac 12
\varepsilon_{10}(m)=\varepsilon_1(M),\ee then $z=0$ is a regular
point of $v$.
\end{pro}
\textsc{Proof} By condition (\ref{332}), we can find a number
$r_1\in ]0,1]$ such that
$$\sup\limits_{0<r\leq r_1}E(r)\leq m$$ and scale $v$ and $p$ so that
$$u(x,t)=r_1v(r_1x,r_1^2t),\qquad q(x,t)=r^2_1p(r_1x,r_1^2t),\qquad |x|<1, \quad
-1<t<0.$$ The pair $u$ and $q$ is then a suitable weak solution to
the NSE's in $Q$ and
$$\sup\limits_{0<r\leq 1}E(u,r)\leq m$$
and \be\la{334}g(u,q)<\frac 12 \varepsilon_{10}(m).\ee By
(\ref{334}), among  quantities $\lim\inf_{r\to 0}E(u,r)$,
$\lim\inf_{r\to 0}A(u,r)$,
 $\lim\inf_{r\to 0}\\C(u,r)$,
$\lim\inf_{r\to 0}H(u,r)$, and $\lim\inf_{r\to 0}D_0(q,r)$, there
should be at least one, which is less than
$1/2\varepsilon_{10}(m)$. For example, let $\lim\inf_{r\to
0}E(u,r)$ do so. Then we can find a number $r_0 \in]0,\min\{1/4,
(A^\frac 32(u,1)+D_0^2(u,1))^{-2}\}[$ such that
$$E(u,r_0)<\frac 34\varepsilon_{10}(m)$$
and thus
$$g_{r_0}(u,q)<
\varepsilon_{10}(m).$$ By Proposition \ref{3p1}, the point $z=0$
is a regular point of $u$ and therefore it is a regular point of
$v$. Proposition \ref{3p2} is proved.

In the same way, one can prove the following statements.
\begin{pro}\la{3p3}Let the pair $v$ and $p$ be a suitable weak solution
to the NSE's in $Q$.

For any $M>0$, there exists a positive number $\varepsilon_2(M)$
such that if
$$\limsup\limits_{r\to 0} A(r)<M$$
and
$$g<\varepsilon_2(M),$$
then $z=0$ is a regular point of $v$.\end{pro}
\begin{pro}\la{3p4}Let the pair $v$ and $p$ be a suitable weak solution
to the NSE's in $Q$.

For any $M>0$, there exists a positive number $\varepsilon_3(M)$
such that if
$$\limsup\limits_{r\to 0} C(r)<M$$
and
$$g<\varepsilon_3(M),$$
then $z=0$ is a regular point of $v$.\end{pro}

\textsc{Proof of Theorem \ref{1t4}} It is a direct consequence of
Propositions \ref{3p2}--\ref{3p4}. Theorem \ref{1t4} is proved.

Let us give some comments on the proof of Theorem \ref{1t7}. The
crucial point is an analog of Proposition \ref{3p1}. It can be
formulated as follows.
\begin{pro}\la{3p5} Let the pair $v$ and $p$ be a suitable weak
solution to the NSE's in $Q$.

For any $M>0$, there exists a positive number
$\widehat{\varepsilon}_{10}=\widehat{ \varepsilon}_{10}(M)$ with
the following property. If
 $$   \sup\limits_{0<r\leq 1} E(r)=E_{0}\leq M$$
and
  $$  E_3(r_*) <\widehat{\varepsilon}_{10}(M)$$
for some $r_*\in ]0,\min\{1/4,(A^\frac 32(1)+D_0^2(1))^{-2}\}[$,
then $z=0$ is a regular point of $v$.\end{pro} \textsc{Proof} We
repeat the proof  of Proposition \ref{3p1} with the following
modifications. Instead of (\ref{34}), we have
$$E_3(v^n,r_n)=\frac 1{r_n}\int\limits_{Q(r_n)}|v^n_{,3}|^2dz\leq
\frac 1n.$$ Now, (\ref{310}) should be replaced with
$$E_3(u^{n_k},1)\to 0\qquad as \quad n\to +\infty.$$
The latter implies
$$u^{n_k}_{,3}\to 0\qquad in \quad L_2(Q)$$
and thus $u_{,3}=0$ in $Q$. This makes it possible to introduce a
new function  $U(x_1,x_2,t)=u(x_1,x_2,x_3,t)$ which satisfies the
inequalities
$$E(u,r)\geq \frac 1r\int\limits_{-r^2}^0dt\int\limits^{r/\sqrt{2}}
_{-r/\sqrt{2}}dx_3\int\limits_{x^2_1+x^2_2<r^2/2}|\na_2
U|^2dx_1dx_2=$$
$$=\sqrt{2}\int\limits_{-r^2}^0\int\limits_{x^2_1+x^2_2<r^2/2}|\na_2
U|^2dx_1dx_2dt\geq \frac 1{\sqrt{2}}E(u,r/\sqrt{2})$$ for any
$0<r\leq 1$. Here, $\na_2$ is the two-dimensional gradient.

So, we can state that
$$\int\limits_{-1}^0\int\limits_{x^2_1+x^2_2<1/2}|\na_2
U|^2dx_1dx_2dt\leq E(u,1)<+\infty$$ and thus
$$\int\limits_{-r^2}^0\int\limits_{x^2_1+x^2_2<r^2/2}|\na_2
U|^2dx_1dx_2dt\to 0\qquad as \quad r\to 0.$$ This means that
$E(r)\to 0$ as $r\to 0$ and, by Theorem \ref{1t2}, $z=0$ is a
regular point of $u$. So, we have shown the validity of
(\ref{323}) and  (\ref{325}).

Next, we proceed in a slightly different way. By (\ref{24}), we
have
$$D_0(q^{n_k},r)\leq c\Big[\Big(\frac r\varrho\Big)^\frac 52D_0(q^{n_k},\varrho)
+\Big(\frac \varrho r\Big)^2C(u^{n_k},\rho)\Big]$$ for all
$0<r<\varrho\leq 1$. Thanks to (\ref{312}),  the new version of
the previous inequality can be given. It has the form
$$D_0(q^{n_k},r)\leq c\widetilde{d}_0(M)\Big[\Big(\frac r\varrho\Big)^\frac 52
+\Big(\frac \varrho r\Big)^2C(u^{n_k},\rho)\Big]$$ for all
$0<r<\varrho\leq 1$. Passing to the limit as $k\to +\infty$, we
find from (\ref{315}) the estimate
$$\limsup\limits_{k\to 0}D_0(q^{n_k},r)\leq c\widetilde{d}_0(M)\Big[\Big(\frac r\varrho\Big)^\frac 52
+\Big(\frac \varrho r\Big)^2C(u,\rho)\Big]$$ for all
$0<r<\varrho\leq 1$. Setting $\varrho=\gamma r\leq r_1$, we derive
from (\ref{323}) the inequality
$$\limsup\limits_{k\to 0}D_0(q^{n_k},r)\leq c\widetilde{d}_0(M)\Big[\ga^{-\frac
52} +c\ga^2d^3_3(M)\ga^3r^3]$$
$$\leq c\widetilde{d}_4(M)\Big[\ga^{-\frac
52} +\ga^5r^3]$$ which is valid for any $\ga>1$ and for any $r>0$
satisfying the condition $\ga r\leq r_1$. So, we have (see
(\ref{325}))
$$cd^3_3(M)r^3+c\widetilde{d}_4(M)\ga^{-\frac
52}+c\widetilde{d}_4(M)\ga^5r^3\geq \varepsilon $$ for the same
$\ga$ and $r$ as in the previous inequality.
Let us choose $\ga_0(\varepsilon,M)>1$ and fix it so that
$$c\widetilde{d}_4(M)\ga_0^{-\frac
52}(\varepsilon,M)\leq \varepsilon/4.$$ Then
$$cd^3_3(M)r^3+c\widetilde{d}_4(M)\ga_0^5(\varepsilon,M)r^3\geq 3\varepsilon/4 $$
for all $0<r\leq r_1/\ga_0(\varepsilon,M)$. The latter inequality
does not hold for sufficiently small $r$. Proposition (\ref{3p5})
is proved.

The remaining part of the proof of Theorem \ref{1t7} is the same
as in  Theorem \ref{1t4}.

 \setcounter{equation}{0}
\section {Appendix }

 Consider the  initial boundary value problem
 for the following system of linear equations
 \be\la{l1}\pa_t v+u\cdot \na v-\De v=f \qquad in \quad Q.\ee
We assume that unknown vector-valued function $v:Q\to\mathbb R^3$
satisfies the homogeneous conditions on the parabolic boundary
$\pa' Q$, i.e.,
\begin{equation}\label{l2}
    v|_{\pa'Q}=0.
\end{equation}
Here,
$u:Q\to\mathbb R^3$ and $f:Q\to\mathbb R^3$ are given functions
satisfying the following conditions
\begin{equation}\label{l3}
    u\in L_2(Q), \qquad \div u=0\qquad in \quad  Q,
\end{equation}
and
\begin{equation}\label{l4}
    f\in L_\frac {5}3(Q).
\end{equation}
\begin{definition}\la{ld1} The function $v$ is called a weak solution
to initial boundary value problem (\ref{l1}) and (\ref{l2}) if it
satisfies the conditions:
\begin{equation}\label{l5}
    v\in L_{2,\infty}(Q)\cap L_2(-1,0;{\stackrel{\circ}{W}}{^1_2}(B))
\end{equation}
\begin{equation}\label{l6}
    t\to\int\limits_{B}v(x,t)\cdot w(x)dx\,\,\mbox{is continuous
    on}\,\,[-1,0]\,\,\mbox{for any}\,\,w\in L_2(B);
\end{equation}
\begin{equation}\label{l7}
    \int\limits_Q\Big(-v\cdot\pa_t w-v\otimes u:\na w+
    \na v:\na w\Big)dz= \int\limits_Qf\cdot wdz
\end{equation}
for any $w\in C^\infty_0(Q)$;
\begin{equation}\label{l8}
    \|v(\cdot,t)\|_{L_2(B)}\to 0
\end{equation}
as $t\to -1+0$;
\begin{equation}\label{l9}
    \frac 12\int\limits_{B}|v(x,t)|^2dx+
    \int\limits_{-1}^t\int\limits_{B}|\na
    v|^2dxdt'\leq \int\limits_{-1}^t\int\limits_{B}f\cdot vdxdt'
\end{equation}
for all $t\in [-1,0]$.
\end{definition}
\begin{theorem}\la{lt2} Assume that conditions (\ref{l3}) and (\ref{l4})
hold. There exists at least one weak solution $v$ to initial
boundary value problem (\ref{l1}), (\ref{l2}). In addition, it has
the following differentiability properties
\begin{equation}\label{l10}
v\in L_{3,\infty}(Q)\cap L_5(Q).
\end{equation}
Moreover, if
\begin{equation}\label{l11}
    u\in L_{\frac {10}3}(Q),\end{equation}
    then problem (\ref{l1}), (\ref{l2}) has a unique weak
    solution.
\end{theorem}
\textsc{Proof} We start with the proof of the existence. It is
easy to find smooth functions $u^\de$ and $f^\de$ with the
following properties
\begin{equation}\label{l12}
    \div u^\de=0,  \qquad u^\de\to u\quad in \quad L_2(Q)
\end{equation}
and
\begin{equation}\label{l13}
    f^\de\to f\quad in \quad L_{\frac {5}3}(Q)
\end{equation}
as $\de\to 0$. As it is shown in \ci{LSU}, see also \ci{KL}, there
is a unique smooth solution $v^\de$ to the initial boundary value
problem:
 \be\la{l14}\pa_t v^\de+u^\de\cdot \na v^\de-\De v^\de=f^\de \qquad in \quad Q\ee
and
\begin{equation}\label{l15}
    v^\de|_{\pa'Q}=0.
\end{equation}
\begin{lemma}\la{ll3} Assume that we are given three
sufficiently smooth functions $v$, $u$, and $f$ satisfying
(\ref{l1})
 and (\ref{l2}).
 Then, for all $t\in ]-1,0[$, the
following inequalities are valid:
$$\pa_t\int\limits_B|v|^2dx+\int\limits_B|\na v|^2dx
\leq c\Big[\Big(\int\limits_B|f|^\frac 53dx\Big)^\frac 12+
\Big(\int\limits_B|f|^\frac 53dx$$
\begin{equation}\label{l16}
+\int\limits_B|\div u|^\frac 52dx\Big)\int\limits_B|v|^2dx\Big]
\end{equation}
and $$ \pa_t\int\limits_B|w|^2dx+\int\limits_B|\na w|^2dx\leq
c\Big[ \int\limits_B|f|^\frac 53dx
\Big(\int\limits_B|w|^2dx\Big)^\frac
4{9}$$\begin{equation}\label{l17}+\int\limits_B|\div u|^\frac 52dx
\int\limits_B|w|^2dx\Big],
\end{equation}
where $w=|v|^\frac 32$.\end{lemma}
 \textsc{Proof} Multiplying identity
(\ref{l1}) by $v$ and integrating by parts, we have
\be\la{l18}\frac 12\pa_t\int\limits_B|v|^2dx+\int\limits_B|\na
v|^2dx= \int\limits_Bf\cdot vdx+\frac 12\int\limits_B|v|^2\div udx
=I_1+I_2.\ee Now, our aim is to estimates $I_1$ and $I_2$. For
$I_1$, we use H\"older's inequality
$$I_1\leq \Big(\int\limits_B|f|^\frac 53dx\Big)^\frac 35
\Big(\int\limits_B|v|^\frac 52dx\Big)^\frac 25.$$ The second
multiplier on the right hand side of the latter relation can be
estimated with the help of the multiplicative inequality
$$I_1\leq c\Big(\int\limits_B|f|^\frac 53dx\Big)^\frac 35
\Big(\int\limits_B|v|^2dx\Big)^\frac 7{20}
 \Big(\int\limits_B|\na v|^2dx\Big)^\frac 3{20}.$$
Next, by Young's inequality,
$$I_1\leq \frac 14\int\limits_B|\na v|^2dx+
c\Big(\int\limits_B|f|^\frac 53dx\Big)^\frac {12}{17}
\Big(\int\limits_B|v|^2dx\Big)^\frac 7{17}.$$ Since $5/17
+7/17=12/17$, one may apply Young's inequality once more. As a
result, we have
\begin{equation}\label{l19}
I_1\leq \frac 14\int\limits_B|\na v|^2dx+c\int\limits_B|f|^\frac
53dx\int\limits_B|v|^2dx+c\Big(\int\limits_B|f|^\frac
53dx\Big)^\frac 12.\end{equation}

Now, we are going to evaluate $I_2$. H\"older's inequality gives
$$I_2\leq \frac 12\Big(\int\limits_B|v|^\frac {10}3
dx\Big)^\frac 35\Big(\int\limits_B|\div u|^\frac 52dx\Big)^\frac
{2}{5}.$$ Then we use  a multiplicative inequality
$$I_2\leq c\Big(\int\limits_B|\na v|^2dx\Big)^\frac 35
\Big(\int\limits_B| v|^2dx\Big)^\frac 25\Big(\int\limits_B|\div
u|^\frac 52dx\Big)^\frac {2}{5}.$$ And, by Young's inequality,
\begin{equation}\label{l20}
    I_2\leq \frac 14\int\limits_B|\na v|^2dx+c\int\limits_B|\div
u|^\frac 52dx\int\limits_B| v|^2dx.
\end{equation}
Now, (\ref{l16}) follows from (\ref{l18})-(\ref{l20}).

To prove the second estimate, we proceed as follows. First, we
multiply identity (\ref{l1}) by $|v|v$ and then integrate by
parts. As a result, we have the identity
$$
    \frac 13\pa_t\int\limits_B| v|^3dx+
    \int\limits_B| v|(|\na v|^2+|\na
    |v||^2)dx=\int\limits_B|v|f\cdot vdx$$$$+\frac 13\int\limits_B|v|^3
    \div u dx,
$$
which implies the inequality
\begin{equation}\label{l21}
 \frac 13\pa_t\int\limits_B| w|^2dx+
 \frac 49\int\limits_B|\na w|^2dx\leq J_1+J_2,
\end{equation}
where
$$J_1=\Big(\int\limits_B| f|^\frac 53dx\Big)^\frac 35
\Big(\int\limits_B| w|^\frac {10}3dx\Big)^\frac 25,\qquad
J_2=\frac 13\int\limits_B|w|^2
    \div u dx$$
The term $J_2$ is estimated in the same way as $I_2$. So,
\begin{equation}\label{l22}
    J_2\leq \frac 1{18}\int\limits_B|\na w|^2dx+c\int\limits_B|\div
u|^\frac 52dx\int\limits_B| w|^2dx.
\end{equation}
For $J_1$, one can apply a multiplicative inequality and Young's
inequality and find
$$J_1\leq c\Big(\int\limits_B| f|^\frac 53dx\Big)^\frac 35
\Big(\int\limits_B| w|^2dx\Big)^\frac 4{15} \Big(\int\limits_B|
\na w|^2dx\Big)^\frac 25$$
$$\leq \frac 1{18}\int\limits_B|\na w|^2dx+\int\limits_B| f|^\frac 53dx
\Big(\int\limits_B| w|^2dx\Big)^\frac 4{9}.$$ Combining the latter
estimate with (\ref{l21}) and (\ref{l22}), we complete the proof
of (\ref{l18}). Lemma \ref{ll3} is proved.

Now, we proceed with the proof of the theorem. From Lemma
\ref{ll3} and from Gronwall's lemma, it follows that certain norms
of $v^\de$ and $w^\de=|v^\de|^\frac 32$ are
 bounded uniformly with respect to $\de$:
\begin{equation}\label{l23}
    \|v^\de\|_{L_{3,\infty}(Q)} +\|w^\de\|_{L_{2,\infty}(Q)}
     +\|\na v^\de\|_{L_2(Q)}+\|\na w^\de\|_{L_2(Q)}\leq c<+\infty.
\end{equation}
Moreover, by a multiplicative inequality and by (\ref{l23}),
\begin{equation}\label{l24}
\|v^\de\|_{L_5(Q)}\leq c<+\infty
\end{equation}
and thus
\begin{equation}\label{l25}
    \||v^\de||u^\de|\|_{L_ {10/7}(Q)}\leq c<+\infty.
\end{equation}
The latter, together with identity (\ref{l7}), allows us to state
that
\begin{equation}\label{l26}
    \|\pa_t v^\de\|_{L_ {10/7}(-1,0;W^{-1}_{ 10/7}(B))}\leq
    c<+\infty.
\end{equation}
Now, the existence of at least one weak solution to initial
boundary value problem (\ref{l1}) and (\ref{l2}) can be deduced
from estimates (\ref{l23})-(\ref{l26}) in a more or less standard
way.

Now, let us switch to the proof of the uniqueness of weak
solutions to initial boundary value problem (\ref{l1}) and
(\ref{l2}) in the class of weak solutions in the sense of
Definition \ref{ld1}. We start with a simple remark
\begin{equation}\label{l27}
  {\stackrel{\circ}{W}}{^1_2}(Q)\cap L_5(Q)=[C^\infty_0(Q)]_
  {W^1_2(Q)\cap L_5(Q)}.
\end{equation}
Assume that $v^1$ is another solution to initial boundary value
problem (\ref{l1}) and (\ref{l2}). By (\ref{l11}),
\be\la{l28}u\cdot \na v^1\in L_\frac 54(Q).\ee Taking into account
(\ref{l27}) and (\ref{l28}), let us transform identity (\ref{l7})
into the form:
\begin{equation}\label{l29}
 \int\limits_Q\Big(-v^1\cdot\pa_t w+ (u\cdot\na v^1)\cdot w+
    \na v:\na w\Big)dz= \int\limits_Qf\cdot wdz
\end{equation}
for any $w\in {\stackrel{\circ}{W}}{^1_2}(Q)\cap L_5(Q)$.

On the other hand, according to the first part of the theorem, we
have
\begin{equation}\label{l30}
    |v||u|\in L_2(Q)
\end{equation}
and thus
\begin{equation}\label{l31}
 \int\limits_Q\Big(-v\cdot\pa_t w-v\otimes u:\na w+
    \na v:\na w\Big)dz= \int\limits_Qf\cdot wdz\end{equation}
for any $w\in  {\stackrel{\circ}{W}}{^1_2}(Q)$.

Fix an arbitrary number $\varepsilon\in ]0, 1/4[$ and take a
nonnegative smooth cut-off function $\chi$ having a support
 in $]-1+\varepsilon/2,-\varepsilon/2[$. Let $\omega_\varrho$ be a
 standard smoothing kernel and let
$$v_\varrho(x,t)=\int\limits_{-1}^0\omega_\varrho(t-s)v(x,s)ds.$$
Now, we introduce two test functions: $$w=(\chi
v_\varrho)_\varrho,\qquad w^1=(\chi v^1_\varrho)_\varrho.$$ For
sufficiently small $\varrho$, we have
\begin{equation}\label{l32}
    w\in {\stackrel{\circ}{W}}{^1_2}(Q)\cap L_5(Q),\qquad
    w^1\in  {\stackrel{\circ}{W}}{^1_2}(Q).
\end{equation}
Next, we use the first test function in (\ref{l30}) and then both
of them
 in (\ref{l31}). Taking into account well-known
properties of smoothing kernel, we can produce three identities:

$$
\int\limits_Q\Big(-v^1_\varrho\cdot \pa_t(\chi v_\varrho)+(u\cdot
\na v^1)_\varrho\cdot \chi v_\varrho+\na v^1_\varrho:\chi\na
v_\varrho
\Big)dz$$\begin{equation}\label{l33}=\int\limits_Qf_\varrho\cdot
v_\varrho dz,
\end{equation}
$$\int\limits_Q\Big(-v_\varrho\cdot \pa_t(\chi v^1_\varrho)-(v\otimes
 u)_\varrho: \chi \na v^1_\varrho+\na v_\varrho:\chi\na
v^1_\varrho \Big)dz$$
\begin{equation}\label{l34}=\int\limits_Qf_\varrho\cdot
v^1_\varrho dz,
\end{equation}
and
$$\int\limits_Q\Big(-v_\varrho\cdot \pa_t(\chi v_\varrho)-(v\otimes
 u)_\varrho: \chi \na v_\varrho+\na v_\varrho:\chi\na
v_\varrho \Big)dz$$
\begin{equation}\label{l35}=\int\limits_Qf_\varrho\cdot
v_\varrho dz.
\end{equation}

Add (\ref{l33}) and (\ref{l34}) and integrate the sum by parts, we
find
$$\int\limits_Q\Big(-v_\varrho\cdot v^1_\varrho\pa_t\chi+(u\cdot
\na v^1)_\varrho\cdot \chi v_\varrho-(v\otimes
 u)_\varrho: \chi \na v^1_\varrho $$$$+2\chi\na v_\varrho:\na
v^1_\varrho
\Big)dz=\int\limits_Qf_\varrho\cdot(v_\varrho+v^1_\varrho)dz.$$
Passing to the limit as $\varrho\to 0$, we arrive at the identity
\begin{equation}\label{l36}
\int\limits_Q\Big(-v\cdot v^1\pa_t\chi+2\chi\na v:\na v^1 \Big)dz
=\int\limits_Qf\cdot(v+v^1)dz.
\end{equation}
Now, we argue as follows. Taking into account two properties of
weak solutions (\ref{l6}) and (\ref{l8}) and choosing our cut-off
function $\chi$ in an appropriate way, we can state that
\be\la{l37}\int\limits_Bv(x,t)\cdot v^1(x,t)dx+2\int\limits_{-1}^t
\int\limits_B\na v:\na v^1dxdt'=\int\limits_{-1}^t \int\limits_B
f\cdot(v+v^1)dxdt'\ee for a.a. $t\in ]-1,0[$. Proceeding in the
same way, we derive from (\ref{l35}) another identity
\begin{equation}\label{l38}
\frac 12\int\limits_B|v(x,t)|^2dx+\int\limits_{-1}^t
\int\limits_B|\na v|^2dxdt'=\int\limits_{-1}^t \int\limits_B
f\cdot vdxdt'
\end{equation}
for a.a. $t\in ]-1,0[$. But $v^1$ is a weak solution and therefore
it obeys global energy inequality (\ref{l9})
\begin{equation}\label{l39}
    \frac 12\int\limits_{B}|v^1(x,t)|^2dx+
    \int\limits_{-1}^t\int\limits_{B}|\na
    v^1|^2dxdt'\leq \int\limits_{-1}^t\int\limits_{B}f\cdot v^1dxdt'
\end{equation}
for all $t\in [-1,0]$.  From (\ref{l37})--(\ref{l39}), it is easy
to derive the estimate
 $$\frac 12\int\limits_{B}|v^1(x,t)-v(x,t)|^2dx+
    \int\limits_{-1}^t\int\limits_{B}|\na
    (v^1-v)|^2dxdt'\leq 0$$
for a.a. $t\in ]-1,0[$. The latter implies $v^1=v$. Theorem
(\ref{lt2}) is proved.

Now, let us consider the pair $v$ and $p$ being a suitable weak
solution to the Navier-Stokes equations in $Q$, see definition
\ref{1d1}. Our aim is


\begin{pro}\la{ap1}Assume that the pair $v$ and $p$ is a suitable
weak solution in $Q$. Suppose that $p$ is independent of spatial
variables $x$. Then $v$ is H\"older continuous in the
 $\overline{Q}(3/8)$.\end{pro}
\textsc{Proof} Fix a smooth cut-off function $\psi$ vanishing in a
neighborhood of $\pa'Q$ and being equal to one in
$\overline{Q}(3/4)$. Let
$$\overline{v}=\psi v, \qquad f=v\pa_t \psi-2\na v\na \psi-v\De \psi+
v\cdot \na \psi v, \qquad u=v.$$ According to the theory of
multiplicative inequalities,
\begin{equation}\label{a4}
    u\in L_\frac {10}3(Q),\qquad \div u=0.
\end{equation}
Moreover,
\begin{equation}\label{a5}
    \overline{v}\in L_{2,\infty}(Q)\cap L_2(-1,0;{\stackrel{\circ}{W}}{^1_2}(B))
\end{equation}
and
\begin{equation}\label{a6}
 \int\limits_Q\Big(-\overline{v}\cdot\pa_t w-\overline{v}\otimes u:\na w+
    \na \overline{v}:\na w\Big)dz= \int\limits_Qf\cdot wdz
\end{equation}
for any $w\in C^\infty_0(Q)$.

Now, let us prove that $\overline{v}$ is a weak solution to the
initial boundary value problem (\ref{l1}) and (\ref{l2}). By
(\ref{a5}) and (\ref{a6}), it remains to verify that
$\overline{v}$ satisfies (\ref{l6}) and (\ref{l9}). Obviously,
(\ref{l8}) is fulfilled since $\overline{v}$ is identically equal
to zero in a neighborhood of $t=0$.

To check (\ref{l6}), we note that
$$\overline{v}\otimes u\in L_\frac 53(Q).$$
This, together with (\ref{a6}), gives us the estimate for
$$\|\pa_t\overline{v}\|_{L_{5/3}(-1,0;W^{-1}_{5/3}(B))}<+\infty.$$
Then, (\ref{l6}) follows from (\ref{a5}).

Next, integrating by parts, we can derive the following identity
$$I(t)=\int\limits_B|\overline{v}(x,t)|^2dx+2\int\limits_{-1}^t
\int\limits_B|\na \overline{v}|^2dxdt'-2\int\limits_{-1}^t
\int\limits_B f\cdot\overline{v}dxdt'$$
$$=\int\limits_B\psi^2(x,t)|{v}(x,t)|^2dx+2\int\limits_{-1}^t
\int\limits_B\psi^2|\na{v}|^2dxdt'$$$$ -\int\limits_{-1}^t
\int\limits_B\Big(|v|^2(\De \psi^2+\pa_t\psi^2)+v\cdot\na
\psi^2|v|^2 \Big)dxdt'.$$ By assumptions of Proposition \ref{ap1},
$I(t)\leq 0$ for a.a. $t\in ]-1,0[$. Here, we have used the
identity
$$\int\limits_Bv\cdot \na \psi^2pdx=0.$$
However, (\ref{l6}) means that the  inequality $I(t)\leq 0$ holds
for all $t\in [-1,0]$. So, we have demonstrated that
$\overline{v}$ is a weak solution to initial boundary value
problem (\ref{l1}) and (\ref{l2}). Then, according to Theorem
\ref{lt2}, we can state that
$$\overline{v}\in L_{3,\infty}(Q)$$ and thus \be\la{a7} v\in
L_{3,\infty}(Q(3/4)).\ee In this case, we have a suitable weak
solution $v$ and $p$ in $Q(3/4)$ which satisfies additional
condition (\ref{a7}). Using scaling and results of \ci{ESS4}, we
show that $v$ is H\"older continuous in the closure of the set
$Q(3/8)$. Proposition \ref{ap1} is proved.

G. Seregin\\
Steklov Institute of Mathematics at St.Petersburg, \\
Fontanka 27, 191023 St.Peterburg, Russia,\\
seregin@pdmi.ras.ru

\end{document}